\documentclass[reqno, 12pt]{amsart}

\usepackage{amssymb}
\usepackage{amsfonts}
\usepackage{epsfig}
\usepackage{color}
\usepackage{epsf}
\usepackage{float}
\usepackage{a4wide}

\newtheorem{thm}{Theorem}
\newtheorem{lemma}[thm]{Lemma}

\newtheorem{propos}[thm]{Proposition}

\def\reff#1{(\ref{#1})}
\def\one{{\mathbf 1}} \def\cX{{\mathcal X}}\def\cM{{\mathcal M}} \def\cL{{\mathcal L}}
 
\def\ualpha{{\underline\alpha}} \def\ueta{{\underline\eta}}
\def\tmu{{\tilde\mu}}
\def\tnu{{\tilde\nu}}
\def\teta{{\tilde\eta}}
\def\vep{{\varepsilon}}
\def\uell{{\underline\ell}}
\def\had{{\textsc{had}}}
\def\tasep{{\textsc{tasep}}}

\def\R{{\mathbb R}}  
\def\P{{\mathbb P}}  
\def\Z{{\mathbb Z}}  
\def\E{{\mathbb E}}  

\def\square{\ifmmode\sqr\else{$\sqr$}\fi}
\def\sqr{\vcenter{
         \hrule height.1mm
         \hbox{\vrule width.1mm height2.2mm\kern2.18mm\vrule width.1mm}
         \hrule height.1mm}}                  
\def\Xnup{\cX^{n\uparrow}}

\begin{document}



\centerline{\Large Multiclass Hammersley-Aldous-Diaconis process}
\centerline{\Large and multiclass-customer queues}

\vskip 3mm

\centerline{\bf Pablo A. Ferrari}
\centerline{\it Universidade de S\~ao Paulo}
\vskip 2mm
\centerline{\bf James B. Martin}
\centerline{\it University of Oxford}


\vskip1truecm\rm
\noindent {\bf Abstract:} 
In the Hammersley-Aldous-Diaconis process infinitely many particles sit in $\R$
and at most one particle is allowed at each position. A particle at $x$, whose 
nearest neighbor to the right is at $y$, jumps at rate $y-x$ to a position
uniformly distributed in the interval $(x,y)$. The basic coupling between
trajectories with different initial configuration induces a process with
different classes of particles. We show that the invariant measures for the
two-class process can be obtained as follows.  First, a stationary $M/M/1$ queue
is constructed as a function of two homogeneous Poisson processes, the arrivals
with rate $\lambda$ and the (attempted) services with rate $\rho>\lambda$. Then
put the first class particles at the instants of departures (effective services)
and second class particles at the instants of unused services. The procedure is
generalized for the $n$-class case by using $n-1$ queues in tandem with $n-1$
priority-types of customers. A \emph{multi-line} process is introduced; it
consists of a coupling (different from Liggett's basic coupling), having as
invariant measure the product of Poisson processes. The definition of the
multi-line process involves the \emph{dual points} of the space-time Poisson
process used in the graphical construction of the system. The coupled process is a
transformation of the multi-line process and its invariant measure the
transformation described above of the product measure.


\paragraph{\bf AMS-Classification} 60K35 60K25 90B22

\section{The Hammersley-Aldous-Diaconis process}
\label{i}
The state space $\cX$ is an appropriate subset of locally finite subsets of
$\R$. Elements of $\cX$ are called configurations and elements of a configuration
$\eta\in\cX$ are called \emph{particles}. The \emph{Hammersley-Aldous-Diaconis
  (\had{} ) process} (Hammersley \cite{H72}, Aldous and Diaconis \cite{AD95}) can be
informally described by saying that a particle sitting at $r\in\R$ waits an
exponentially distributed random time with rate equal to the distance to the
nearest particle to its right, located at $r'>r$ (say) to jump to a site
uniformly distributed in the interval $[r,r']$.  Alternatively, bells ring at
space-time points at rate 1 and when a bell rings at $(r,t)$, the nearest
particle to the left of $r$ at time $t-$ jumps to $r$.

A Harris graphical construction of the process is the following. Let
$\omega\in\Omega$ be a homogeneous rate-1 Poisson process in the space-time
space $\R\times\R^+$ (or later in $\R^2$), where $\Omega$ is the set of locally
finite subsets of $\R^2$.  We shall use \emph{points} to refer to space-time
elements of $\omega$. If $\eta$ is the particle configuration at time $t-$ and
$(r,t)\in \omega$, then at time $t$ the configuration jumps to $\eta\setminus
\{u\}\cup\{r\}$, where $u=u(\eta,r)$ is the nearest particle in $\eta$ to the
left of $r$.  This construction is well defined in a finite region as the points
can be well ordered by time \cite{AD95}.  The construction in $\R$ was performed
by Aldous and Diaconis \cite{AD95} and then by Sepp\"al\"ainen
\cite{Seppalainen96} in the state space
$\cX=\{\eta:\, \lim_{s\to\infty} |\eta\cup [0,s]|^2/s = \infty\}$ that we adopt.
Here $|\cdot|$ counts the number of elements of a finite set.
Homogeneous Poisson processes in $\R$ give mass 1
to $\cX$.

For fixed initial configuration $\eta$ and points $\omega$, the process
$(\eta_t,\,t\ge 0)$ is a deterministic function of $\eta$ and $\omega$ denoted
$\eta_t=\Phi(t,\eta,\omega)$, $t\ge 0$. In this case we say that the process is
\emph{governed by} $\omega$ with initial configuration $\eta$. It satisfies
  \begin{equation}
    \label{f64}
    \eta_t =
    \Phi(t-s,\eta_s,\tau_s\omega),
\end{equation}
for all $0\le s<t$, where $\tau$ is the time translation operator defined by
$\tau_s\omega=\{(x,t-s),\,(x,t)\in\omega\}$. 
\begin{figure}[htb]
\begin{center}
\input{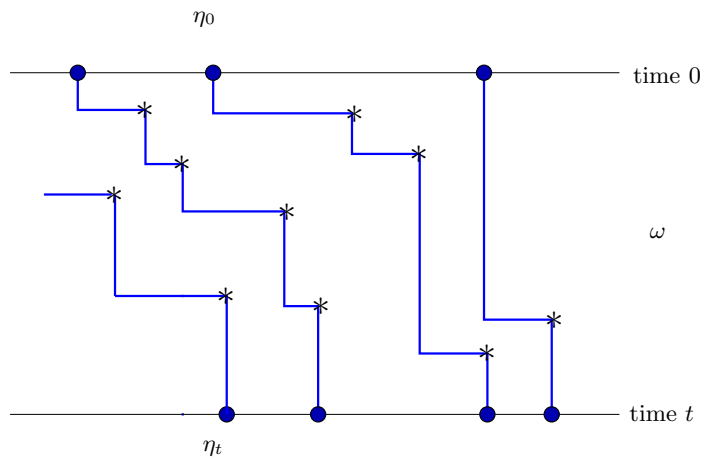}  
\caption{Harris construction. Stars represent the space-time Poisson points in
  $\omega$. Balls represent particle positions at times $0$ and $t$.  }
\label{206}
\end{center}
\end{figure}

\paragraph{\bf Coupled and multiclass process} A joint construction of \had{} 
processes $((\eta^i_t,\, t\ge0),\, i=1,\dots,n)$ with initial configurations
$\eta^1,\dots,\eta^n$ governed by the same points $\omega$ is called
\emph{coupled process}. It is defined by
\begin{equation}
  \label{f1}
  \eta^i_t=\Phi(t,\eta^i,\omega)
\end{equation}
We slightly abuse notation writing $\eta=(\eta^1,\dots,\eta^n)$ and $\eta_t=
\Phi(t,\eta,\omega)\in\cX^n$. In the coupled process when a point is at $(x,t)$,
the closest particle to the left of $x$ in each marginal jumps to $x$
simultaneously.

Consider initial ordered particle configurations
$\eta^1\subset\dots\subset\eta^n$ and run the coupled process. The order is
maintained at later times. Define $\Xnup= \{(\eta^1,\dots,\eta^n)\in\cX:\,
\eta^1\subset\dots\subset\eta^n\}$. Define $R:\Xnup\to\cX^n$ by
\begin{equation} 
  \label{f73}
  (R\eta)^k = \eta^k\setminus\eta^{k-1}
\end{equation}
The process $ \xi_t= R\eta_t$ is called the \emph{multiclass process}.
Particles in $\xi^i_t$ are called $i$-class particles.  The multiclass process
is just a convention to describe a coupled process with ordered initial
configurations. The map $R$ is invertible; its inverse is given by
$(R^{-1}\xi)^k=\xi^1\cup\dots\cup\xi^k$. The process $\xi_t$ governed by
$\omega$ with initial configuration $\xi$ is defined by
\begin{equation}
  \label{f6a}
  \xi_t=\Upsilon(t,\xi,\omega) = R\Phi(t,R^{-1}\xi,\omega).
\end{equation}

\paragraph{\bf Invariant measures}
The Poisson process with density $\lambda$ is an invariant measure for
the \had{}  process for all $\lambda>0$ \cite{AD95}. Our main result is to
construct invariant measures for the multiclass process. The
resulting measure coincides with the law of the departure process of a
stationary multiclass-customer queue system.

Let $A$ and $S$ be particle configurations in $\cX$. Think of $\R$ as time and
construct a continuous time random walk $(Z_r,\,r\in\R)$ jumping one unit up at
times in $A$ and one unit down at times in $S$. We fix $Z_0=0$ but since we are
interested in the increments, the position at a given time is not important. The
increments satisfy
\begin{equation}
  \label{f60}
  Z_r-Z_s= |A \cap   [r,s)|-|S\cap [r,s)|
\end{equation}
Assume that $Z_r\to\mp\infty$ as $r\to\pm\infty$; this implies $Z_r$ visits each
site a finite number of times.  Let $U(A,S)$ be the times in $S$ that $Z_r$
attains a new record down and $D(A,S)$ its complement:
\begin{eqnarray}
  \label{f9}
  U(A,S)&:=& \{r\in S\,:\, Z_r < \inf_{r'<r}
  Z_{r'}\}\nonumber\\
  D(A,S)&:=&S\setminus U(A,S)
\end{eqnarray}

\begin{figure}[htb]
\begin{center}
\input{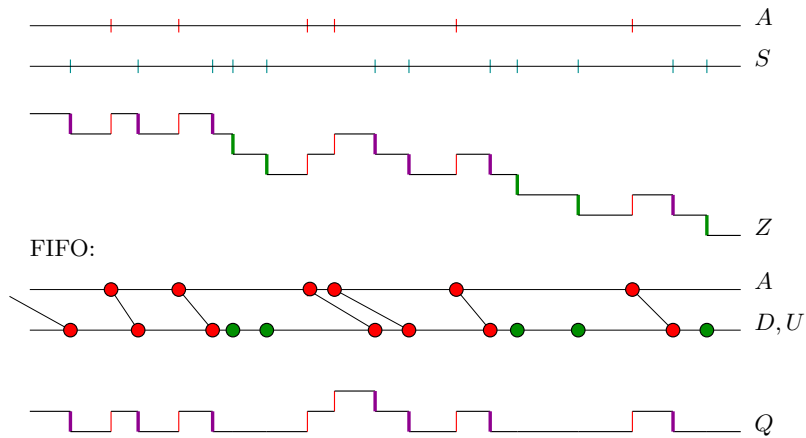}  
\caption{Arrivals and services. Each arrival is linked to its departure time,
  running the queue according to the FIFO (first-in-first-out) schedule.}
\label{145}
\end{center}
\end{figure}

We interpret $A$ as arrival times of customers to a one-server queue and $S$
as the times at which service occurs. At these service times there are
two possibilities: either there is a customer in the system, producing a
\emph{departure}, or there is no customer, in this case there is an
\emph{unused} service.  We collect the departure times in $D$ and the unused
service times in $U$.  If $A$ and $S$ are stationary Poisson processes with
rates $\lambda<\lambda'$, respectively, then $Z_r$ satisfies the conditions
above and $D(A,S)$ are the departure times of a stationary $M/M/1$ queue
$(Q_r\,,\,r\in\R)$ defined by $Q_r = (Z_r - \inf_{r'<r} Z_{r'})^+$; here $Q_r$
is the number of customers in the queue at time $r$. This construction of
$(Q_r,\,r\in\R)$ is standard (reference ??). The process is stationary Markov
and satisfies
\begin{equation}
  \label{f77}
  Q_r-Q_{r-}= \one\{r\in A\}-  \one\{r\in D\,:\,Q_{r-}>0\}
\end{equation}

The operator $D$ has the following properties: (i) $D(\alpha^1,\alpha^2)\subset
\alpha^2$. (ii) If $\tilde{\alpha}^1\subset \alpha^1$ then
$D(\tilde{\alpha}^1,\alpha^2)\subset D(\alpha^1,\alpha^2)$. (iii) If $\alpha^1$
and $\alpha^2$ are independent one-dimensional Poisson processes of densities
$\rho^1$ and $\rho^2$ with $\rho^1<\rho^2$, then $D(\alpha^1,\alpha^2)$ also is
a one-dimensional Poisson process with density $\rho^1$. This is just Burke's
Theorem for a $M/M/1$ queue \cite{BaccelliBremaudbook}.

Let $\nu^\lambda$ be the law of a Poisson process of rate $\lambda$. For
$\rho^1<\dots<\rho^n$ define $\nu := \nu^{\rho^1}\times\dots\times\nu^{\rho^n}$.
Let $\alpha=(\alpha^1,\dots,\alpha^n)\in \cX^n$ be a multi-line configuration
with law $\nu$.  Define a sequence of operators $D^{(n)}:\cX^n\mapsto\cX$ as
follows.  Let $D^{(1)}(\alpha^1)=\alpha^1$, and
then recursively for $n\geq2$, let
\begin{equation}
\label{Dndef}
D^{(n)}(\alpha^1,\alpha^2,\dots,\alpha^n)=
D\left(
D^{(n-1)}(\alpha^1,\dots,\alpha^{n-1}),\alpha^n
\right).
\end{equation}
The configuration $D^{(n)}(\alpha^1,\alpha^2,\dots,\alpha^n)$ represents the
departure process from a system of $(n-1)$ queues in tandem. The arrival process
to the first queue is $\alpha^1$. The service process of the $k$th queue is
$\alpha^{k+1}$, for $k=1,\dots,n-1$.  Finally, for $k=2,\dots,n-1$, the arrival
process to the $k$th queue is given by the departure process of the $(k-1)$st
queue.  This is known as a system of $./M/1$ queues in tandem.

Note $D^{(2)}(\alpha^1,\alpha^2)=D(\alpha^1,\alpha^2)$.  By applying (i)-(iii)
above repeatedly, we obtain that $D^{(n)}(\alpha^1,\dots,\alpha^n) \subset
D^{(n-1)}(\alpha^2,\dots,\alpha^n) \subset \dots \subset \alpha^n$ and if
$\alpha^1,\dots,\alpha^n$ are independent one-dimensional Poisson processes of
densities $\rho^1<\dots<\rho^n$, then $D^{(n)}(\alpha^1,\dots,\alpha^n)$ also is
a one-dimensional Poisson process with density $\rho^1$.

Define the configuration $\eta=(\eta^1,\dots,\eta^n)$ by
\begin{equation}
\label{Tdef}
\eta^k=D^{(n-k+1)}(\alpha^k,\alpha^{k+1},\dots,\alpha^n).
\end{equation}
By construction $\eta^k\subset\eta^{k+1}$ for all $k=1,\dots,n-1$ and for each
$k$, $\eta^k$ has marginal distribution~$\nu^{\rho^k}$.

\begin{figure}[htb]
\begin{center}
\input{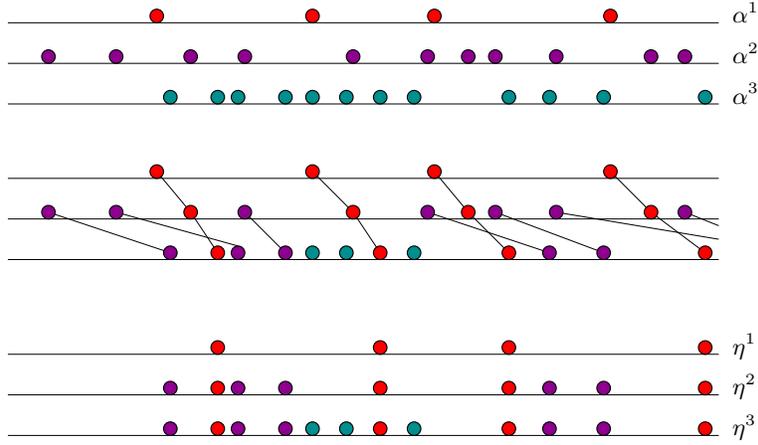} 
\caption{Construction of a coupled configuration $\eta$ and a multiclass
  configuration $\xi$ from a multi-line configuration $\alpha$. $\xi^k$ particles
  are colored according to its class. In the middle picture, arrivals have been
  linked with their departure times in each server, according to the FIFO
  schedule and the customer class. }
\label{176} 
\end{center}
\end{figure}

Define the map $C:\cX^n\mapsto\Xnup$ by $C\alpha=\eta$. Define the measure $\pi$
on $\Xnup$ as the law of $\eta$ (that is, $\pi=C\nu$). Define the multiclass
measure $\mu$ as the law of $\xi=R\eta$, that is, $\mu=RC\nu$. Call $M=RC$. See
Figure \ref{176}.

\begin{thm}
  \label{f5}
  If $\alpha$ has law $\nu$ with $\rho^1<\dots<\rho^n$, then the law $M\nu$ of
  $\xi = M\alpha$ is invariant for the multiclass process $\xi_t$ defined in
  \reff{f6a}.  Equivalently, the law $C\nu$ of $C\alpha$ is invariant for the
  coupled process $\eta_t$ defined in \reff{f1}.
\end{thm}

To prove Theorem \ref{f5} we introduce later a new process
$\alpha_t=(\alpha^1_t,\dots,\alpha^n_t)\in \cX^n$ called \emph{multi-line
  process}. Proposition \ref{f68} shows that the product measure $\nu$ is
invariant for $\alpha_t$ and Proposition \ref{f69} shows that if $\alpha$ is the
initial configuration for the multi-line process $\alpha_t$ then $(C\alpha_t,\,
t\ge 0)$ is the coupled process with initial configuration $C\alpha$.  These
propositions imply that $C\nu$ is invariant for the multiclass process.

The right density of a configuration $\eta\in\cX$ is the limit as $r$ goes to
infinity of the number of $\eta$-particles in $[0,r]$ divided by $r$, if this
limit exists. Analogously define the left density and call simply the
\emph{density} of $\eta$ the result when the left and right density are the
same. The density of an ergodic measure $\nu$ in $\cX$ is the unique value $\lambda$ such
that the set of configurations with density $\lambda$ has $\nu$ probability one.

\paragraph{\bf Uniqueness}
The next result says that the measure $C\nu$ is the unique invariant measure for
the coupled process in the set of ergodic measures on $\cX^n$ with marginal
densities $\rho$.  Its domain of attraction includes the ergodic measures with
densities $\rho$.  Let $\rho^1<\dots<\rho^n$, $\rho=(\rho^1,\dots,\rho^n)$ and
$\rho'=(\rho^1,\rho^2-\rho^1\dots,\rho^n-\rho^{n-1})$. Let $\cM^n(\rho)$ be the
set of ergodic measures on $\cX^n$ such that the $k$th marginal has density
$\rho^k$ for $k=1,\dots,n$.

\begin{thm}
  \label{f5a}
Assume the conditions of  Theorem \ref{f5}. Then,

1) The measure $M\nu$ is the unique invariant measure for the multiclass process
in $\cM^n(\rho')$. The multiclass process starting with a measure in
$\cM^n(\rho')$ converges weakly to $M\nu$ as $t\to\infty$.

2) The measure $C\nu$ is the unique invariant measure for the coupled process in
$\cM^n(\rho)$. The coupled process starting with a measure in $\cM^n(\rho)$
converges weakly to $C\nu$ as $t\to\infty$.
\end{thm}

For all $\lambda>0$ and almost all realizations $\omega$ of the points, 
there is a unique stationary realization of the \had{}  
process at density $\lambda$ governed by
$\omega$. This is shown by the next theorem; notice that the notation $\eta_t$
in this theorem is used for the \had{}  process on $\cX$, that is, with $n=1$.
\begin{thm}
  \label{t3}
  For each $\lambda>0$ there exists an essentially unique function $H_\lambda$
  mapping elements $\omega$ of $\Omega$ to \had{}  trajectories $(\eta_t,t\in\R)$
  such that:
\begin{itemize}
\item[(i)]
The induced law of $(\eta_t,t\in\R)=H_\lambda(\omega)$ is stationary in time.
\item[(ii)]
The marginal law of $\eta_t$ for each $t$ is space-ergodic with particle density
$\lambda$.
\item[(iii)]
With probability 1, $(\eta_t,t\in\R)$ is a \had{}  evolution governed by $\omega$.
\end{itemize}
(Here ``essentially unique'' means that if $H'_\lambda$ is another 
function satisfying the three conditions, then $H_\lambda(\omega)=H'_\lambda(\omega)$
with probability 1).
Then in fact the marginal law of $\eta_t$ for each $t$ is $\nu^\lambda$.
\end{thm}
To stress the dependence on $\lambda$ call $\eta^\lambda_t$ the process
constructed in Theorem \ref{t3}. The construction implies $\eta^\lambda_t\subset
\eta^{\lambda'}_t$ if $\lambda<\lambda'$.  The union of $\eta^\lambda_t$ in
$\lambda$ is a countable dense set of $\R$. It consists on the space coordinates
of the points in $\omega$ with time coordinate less than $t$.

Theorems \ref{f5a} and \ref{t3} are proven in Section \ref{s5}, based on
Proposition \ref{t18} which considers the coupled process $(\eta^1_t,\eta^2_t)$
starting with two independent configurations with the same density $\lambda$. If
one of the configurations is a Poisson process and the other comes from an
ergodic distribution, then, with probability one, the density $c(t)$ of
positions where the two configurations differ at time $t$ is deterministic and
converges to zero as $t$ grows. The proof of the proposition follows an argument
of Ekhaus and Gray \cite{eg93} as implemented by Mountford and Prabhakar
\cite{mp95}.

\paragraph{\bf Multi-line process}
The stationary realization $(\eta^\lambda_t,\,t\in\R)$ of the \had{}  process
governed by $\omega$ and density $\lambda$ of Theorem \ref{t3} induces a point
configuration on $\R\times\R$ consisting on the space-time positions of the
particles just before jumps. Cator and Groeneboom \cite{CG} call them \emph{dual
  points} and prove that they have the same Poisson law as $\omega$ for all
$\lambda>0$. The proof is based on two facts: (a) the reverse \had{}  process with
respect to $\nu^\lambda$ is also a \had{} process, with jumps to the left instead 
of to the right and (b) a
trajectory of the process uniquely determines the points governing it. The dual
points and a new density produce a new stationary realization of the \had{} 
process, which in turn produces new dual points.  Repeating the procedure $n-1$
times, we construct a stationary realization of the \emph{multi-line process}
with $n$ lines.  Fix densities $\rho^1,\dots,\rho^n$ and a realization of the
points $\omega$.  Call $(\alpha^n_t,\,t\in\R)$ the stationary realization of the
\had{}  process at density $\rho^n$ governed by $\omega$ 
(which is well-defined by Theorem \ref{t3}),
and, for $k=n-1,\dots,1$ set
$\omega^{k}=$ dual points generated by $\rho^k$ and $\omega^{k+1}$ with
$\omega^{n+1}=\omega$ (it is convenient for later reasons to start at the
bottom) and $\alpha^{k-1}_\cdot$ as the stationary realization of the \had{} 
process at density $\rho^{k-1}$ governed by $\omega^{k}$.  Define $(\alpha_t,
\,t\in\R)$ as the \emph{multi-line (stationary) process} given by $\alpha_t =
(\alpha^1_t,\dots,\alpha^n_t)$.
\begin{figure}[htb]
\begin{center}
\input{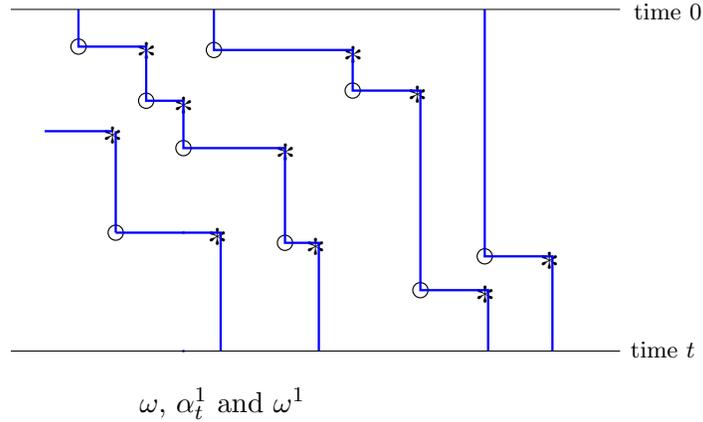}  
\label{299}
\caption{The dual points are represented by circles}
\end{center}
\end{figure}
The generator of the multi-line process is given
later in \reff{f23}. The multi-line process can be constructed for times in
$[0,\infty)$ for any initial distribution; we do so in Section \ref{s2}.  In
this case the $n-1$ first marginals do not follow the \had{}  dynamics unless the
initial distribution is a product of Poisson processes.  We prove two results in
Section \ref{s2}:
\begin{propos}
  \label{f68}
  Let $\rho=(\rho^1,\dots,\rho^k)$ be positive densities. The product
  of Poisson processes $\nu^{\rho^1}\times\dots\times\nu^{\rho^k}$ is
  the unique invariant measure for the multi-line process with marginal
  densities $\rho$.
\end{propos}
 
\begin{propos}
  \label{f69}
  Take a multi-line configuration $\alpha$ such that the multiclass configuration
  $M\alpha$ is well defined. Let $\alpha_t$ be the multi-line process with
  initial configuration $\alpha$ governed by $\omega$. Then $(M\alpha_t,
  \,t\ge0)$ is the multiclass process with initial distribution $M\alpha$ and
  $(C\alpha_t, \,t\ge0)$ is the coupled process with initial distribution
  $C\alpha$, both governed by $\omega$.
\end{propos}
The proof of Theorem \ref{f5} is a consequence of Proposition \ref{f68} (without
the uniqueness part) and Proposition \ref{f69}.  

Pr\"ahofer and Spohn \cite{PS02} and Patrik Ferrari \cite{Patrick04} introduce a
\emph{multilayer dynamics} based on the dual points and related to our multi-line
process by a rotation of $45^o$ (roughly speaking). See the pictures in
\cite{Patrick04}.

\paragraph{\bf Multiclass Burke's Theorem}
The invariant measure for the multiclass \had{}  is also a fixed point for a
multiclass $./M/1$ queue. 

To make this more precise 
we need to introduce more notation in the construction of $\eta$ as
a function $C$ of $\alpha$. As in \cite{FM05} we label multiclass configurations
$\xi^{k,i}$ as a function of the multi-line configuration $\alpha$.  Let
$\xi^{1,1}= \alpha^1$ and define
\begin{eqnarray}
  \label{f2}
  \xi^{k,i} &:=&
  D(\xi^{k-1,i},\alpha^{k}\setminus(\xi^{k,1}\cup\dots\cup
  \xi^{k,i-1}))\nonumber\\ 
  \xi^{k,k}&:=& U(\alpha^{k-1}, \alpha^{k}) 
\end{eqnarray}
for $k=1,\dots,n$ and $i=1,\dots,k$.  Interpret these configurations as the
service times of customers of $k$ classes in the $k$th system; at times in
$\xi^{k-1,i}$ customers of class $i$ arrive to the $k$th system and are served
at times in $\xi^{k,i}$.  The service schedule respects the classes: first class
customers are served first, second class customers are served when there are no
first class customers waiting, etc.  End of services occur at times in
$\alpha^{k}$.  The unused service times of the $k$th system are in $\xi^{k,k}$.
One can think that the $k$th system has infinitely many customers of class $k$
waiting, so that they are served at the service times unused by the lower class
clients.  With this interpretation $\xi^{1,1}=\alpha^1$ are the service times of
a system with infinitely many customers in queue at all times.
Assuming FIFO schedule (first-in-first-out) for customers of the same
class we can link each time in $\alpha^k$ with the corresponding service time in
$\alpha^{k+1}$. This is the meaning of the links in Figure~\ref{176}. In each
line $k$, the particles linked to some particle in the first line will be first
class particles; particles linked to particles in the second line but not to
particles in the first line are second class particles and so on.
 
Let $M:\widetilde \cX^n\to\cX^n$ be the map between multi-line configuration
$\alpha$ and the multiclass configuration $\xi=(\xi^{n,1},\dots,\xi^{n,n})$
using \reff{f9} and \reff{f2}:
\begin{equation}
  \label{f61}
  M\alpha := (\xi^{n,1},\dots,\xi^{n,n})
\end{equation}
where $\widetilde \cX^n\subset \cX^n$ is the subspace of $\cX$ where $M$ is well
defined. This coincides with the map $M=RC$ defined just after \reff{Tdef}. 

The $n$th multiclass stationary queuing system $(Q^n_r,\, r\in\R)$ in
$\{0,1,\dots\}^{\{1,\dots,n-1\}}$ is defined as a deterministic function of the
arrivals $(\xi^{n-1,1},\dots,\xi^{n-1,n-1})$ and the services $\alpha^n$ by
\begin{eqnarray}
  \label{t2}
  Q^n_r(j)-Q^n_{r-}(j) &=& \one\{r\in\xi^{n-1,j}\}-
  \one\{r\in \alpha^n\,:\,Q^n_{r-}(i)=0 \hbox{ for }i<j;\,
  Q^n_{r-}(j)>0\}
\end{eqnarray}
where $Q^n_r(j)$ is the number of customers of class $j$ in the $n$th system at
time $r$. The existence of a process satisfying \reff{t2} can be proved by
induction on $n$ by observing that at times $r$ in $\xi^{n,n}$ the queue is
empty: $Q^n_r(j)=0$ for all $j<n$. The process $Q^n_r$ is not Markov unless the
arrival process is product of homogeneous Poisson and the service process is
also Poisson.  The departure times of class $j$ in system $n$ is given
(equivalently to \reff{f2}) by
\begin{equation}
  \label{t4}
  \xi^{n,j} = \{r\in \alpha^n\,:\,Q^n_{r-}(i)=0\hbox{ for } 
  i<j \hbox{ and } Q^n_{r-}(j)>0\}
\end{equation}
for $j=1,\dots,n-1$. As before $\xi^{n,n}$ is just the set of unused
service times in the $n$th system.

Burke theorem \cite{BaccelliBremaudbook} says that the departures of a
stationary $M/M/1$ queue have the same Poisson law as the arrivals. It applies
to our case for $n=2$: $Q^2_r(1)$ is an $M/M/1$ queue with arrivals $\xi^{1,1}$,
Poisson of rate $\rho^1$ and departures $\xi^{2,1}$. The extension to the
multiclass system says that the departures of the $(n-1)$ classes of the $n$th
system $(\xi^{n,1},\dots,\xi^{n,n-1})$ have the same law as the arrivals
$(\xi^{n-1,1},\dots,\xi^{n-1,n-1})$ to the same system.  This is one of the
consequences of the multiclass Burke's theorem as follows.
\begin{thm}[Multiclass Burke]
  \label{t1} Fix $k\ge 1$ and for $n\ge k$ let $\alpha\in\cX^n$ have law $\nu$
  with $\rho^1<\dots<\rho^n$.  Then the law of $(\xi^{n,1},\dots,\xi^{n,k})\,(=$
  first $k$ coordinates of $M\alpha$) is independent of $n$.
\end{thm}

\paragraph{\bf Proof}
Let $\eta^{(n)}_t$ be the coupled process with initial invariant distribution
$C\nu^{(n)}$, where $\nu^{(n)}$ is the product measure in $\cX^n$ with marginals
$\nu^{\rho^k}$, $\rho^1<\dots<\rho^k$. Since the evolution of the first $k$
coordinates is Markovian, the marginal law of the first $k$ coordinates of
$\eta^{(n)}_t$ under $C\nu^{(n)}$ is invariant for the first $k$ coordinates of
the process $\eta^{(n)}_t$.  By uniqueness of the invariant measure of Theorem
\ref{f5a}, the law of these marginals must coincide with $C\nu^{(k)}$, the
invariant measure for $\eta^{(k)}_t$, the process with $k$ lines and the same
densities. Since $(\xi^{n,1},\dots,\xi^{n,k})=$ first $k$ coordinates of
$R\eta^{(n)}$, its law coincides with the law of $R\eta^{(k)}=M\alpha^{(k)}$,
where $\alpha^{(k)}$ is the multi-line configuration with law $\nu^{(k)}$.
\square

\paragraph{\bf The process in a cycle}
Similar results can be proven for the multiclass process in a cycle $\R_N$. The
space state $\cX_N$ is the set of finite configurations contained in
$[0,N]$. The points $\omega$ are restricted to $\R_N\times\R^+$. When the
configuration at time $t-$ is $\eta_{t-}=\eta$ and there is a point at $(x,t)$,
if there are no $\eta$ particles to the left of $x$, then the rightmost particle
of $\eta$ jumps to $x$. The coupled process is defined as in \reff{f1}.  Take
two finite particle configurations $A_N$ and $S_N$ in $\R_N$ with
$|S_N|>|A_N|$. Extend these configurations periodically to two infinite
configurations $A$ and $S$. Construct $(Z_t,\,t\in \R)$, a periodic process
satisfying \reff{f60}. The resulting periodic configurations of departures
$D(A,S)$ and unused services $U(A,S)$ are periodic and induce configurations
$D_N$ and $U_N$ in $\R_N$.
We construct a (unique)
queue $Q_r$, $r\in \R_N$ such that it satisfies \reff{f77} which has value
$Q_r=0$ for $r\in U_N$; actually $Q_r$ so constructed is the minimal process
satisfying \reff{f77}.

Given a multi-line configuration $\alpha$ in $\cX^n_N$ we construct a multiclass
configuration $\xi=M_N\alpha\in \cX_N^n$, where $M_N$ is defined as the map $M$
substituting $D$ and $U$ by $D_N$ and $U_N$. The analogous to Theorems \ref{f5},
\ref{f5a} and \ref{t1} hold for this process. Letting $\uell
=(\ell^1,\dots,\ell^n)$, the set $\cX^n(\rho)$ of Theorem \ref{f5a} must be
substituted by $\cX^n_N(\uell)$, the set of configurations with exactly $\ell^k$
particles in the $k$th cycle.  The proof of the analogous to Theorem \ref{f5a}
is easy in this case.

\paragraph{\bf Regeneration properties of the invariant measure}
In the two-class invariant measure the second class particles are regeneration
events. More precisely, let $\mu$ on $\cX^2$ be a translation invariant measure
with marginal densities $\rho^1,\rho^2$. The Palm measure ``conditioned to have
a second class particle at the origin'' is the measure $\hat\mu$ defined by
$\hat\mu f= (1/|I|\rho^2)\int d\mu(\xi)\sum_{r\in\xi^2\cap I} f(\theta_r\xi)$,
for any measurable bounded set $I\subset \R$, where $|I|$ is the Lebesgue
measure of $I$ and $\theta_r$ is translation by $r$.  The multiclass invariant
measure $M\nu$ of Theorem \ref{f5} conditioned to have a second class particle
at the origin satisfies: (a) the configuration to the left of the origin is
independent of the configuration to its right; (b) both the distribution of
first class particles to the right of the origin and the positions of the first
plus second class particles to the left of the origin are Poisson processes (or
product measures in the discrete-space case).  These properties were proven by
\cite{DJLStasep} for the two-class invariant measure for the totally asymmetric
simple exclusion process (\tasep); alternative probabilistic proofs can be found
in \cite{FFKtasep}. Angel \cite{Angeltasep} shows that these properties follow
easily from the two-class representation $\xi=M\alpha$ that he called
``collapsing procedure''.  In Section \ref{s7} we show that as in the discrete
case \cite{FM05}, there are \emph{regeneration strings} for the multiclass
invariant measure $M\nu$.

\paragraph{\bf Microscopic shocks}
The process with two classes of particles is a crucial tool to define shock
measures \cite{FKS, Fshock}. Taking a configuration $(\xi^1,\xi^2)$ distributed
according to the invariant measure $M\nu$ for the two-class process conditioned
to have a $\xi^2$ particle at the origin and considering the configuration
consisting on all particles (of any class) to the left of the origin and the
first class particles to the right of it, then the distribution of the resulting
configuration is invariant for the process as seen from an isolated second class
particle.  We discuss this item in Section \ref{s8}.

\paragraph{\bf Multiclass HAD process in $\Z$ and other discrete processes}
The approach has been implemented for a discrete-space multiclass \had{}  process
and other discrete-space processes on $\Z$. In fact, this research started with
the extension of Angel \cite{Angeltasep} result to the multiclass process for
the discrete \had{}  and \tasep{}  we performed in \cite{FM05} and \cite{FMihp}. There
are two main points in this paper: one is to show that the approach works also
in continuous space (the \had{}  process is the most natural model to study \tasep{} 
like questions in $\R$); the other is to show some results that are only
announced in \cite{FMihp}. In particular we show the uniqueness of the
multiclass invariant measure with given marginal densities. The use of dual
points is more intuitive in the continuous \had{}  than in the discrete
processes. Indeed, the dual points in this case are identified by the trajectory
of the process, while in the discrete cases it is necessary to appeal auxiliary
spin-flip processes to identify the dual points.

The the existence of an invariant measure for the two-classes TASEP is first
show by Liggett \cite{Liggettcoupling} \cite{Liggettbook} and then
computed by Derrida, Janowsky and Lebowitz \cite{DJLStasep}, see also
\cite{Speertasep} and \cite{FFKtasep}. Angel \cite{Angeltasep} and Duchi and
Schaeffer \cite{DucSch} description of the two-class invariant measures was the
starting point of our multiclass version \cite{FM05} and \cite{FMihp}.

\section{The multi-line process}
\label{s2}
The multi-line process has configurations
$\alpha_t=(\alpha^1_t,\dots,\alpha^n_t)$ in $\cX^n$. Given a multi-line
configuration $\alpha$ and a position $x\in\R$ let $x^n$ be the position of the
closest particle of $\alpha^n$ to the left of $x$ and inductively for lines
$k=n-1, \dots,1$ let $x^{k}$ be the position of the closest $\alpha^k$ particle
to the left of $x^{k+1}$.  When the configuration at time $t$ is $\alpha$ and
there is an~$\omega$ point in $(x,t)$, the $\alpha^n$-particle at $x^n$ jumps to
$x$ at time $t$ and simultaneously the $\alpha^k$ particle located at $x^{k}$
jumps to $x^{k+1}$ for $k=n-1,\dots,1$.  The function $J:\cX^n\times\R\to\cX^n$
defined by
\begin{equation}
  \label{f22}
  (J(\alpha,x))^k = \alpha^k\cup\{x^{k+1}\}\setminus\{x^{k}\}
\end{equation}
maps the multi-line configuration $\alpha$ before the jumps produced by $x$ to
the configuration after the jumps.  The generator of the process is given by
\begin{equation}
  \label{f23}
  \cL f(\alpha) = \int_\R dx [f(J(\alpha,x))-f(\alpha)]
\end{equation}
This definition is equivalent to the one given in the introduction. For $t$ such
that $(x,t)$ is in $\omega$ denote $x^{n+1}(t)=x$ and $x^k(t)$ the position at
time $t-$ of the particle in the $k$th line jumping due to the point in
$(x,t)$. Then $\omega^k =\{(x^{k}(t),t):\, (x,t)\in\omega\}$ are the dual points
of $\omega^{k+1}$ which in turn govern the process $\alpha^k_t$.

We define another multi-line \had{}  process $\alpha^*_t$ in $\cX^n$ and show that it
is the reverse of $\alpha_t$ with respect to the product of Poisson processes
$\nu$.  For $\alpha^*_t$ the points govern the first line but producing jumps to
the left: when a point appears at $y$, it calls the closest $\alpha^1$ particle
to the right of $y$, located at $y^1$.  Simultaneously the closest $\alpha^2$
particle to the right of $y^1$, located at $y^2$ jumps to $y^1$, and so on.
Calling $y_0=y$, the positions $(y_0,\dots,y_{n})$ are defined as a function of
$\alpha$ and $y$. The multi-line configuration obtained after the jumps produced
by a point at $y$ is called $J^*(\alpha,y)$; its $k$th line is given by
\begin{equation} 
  \label{f24}
  (J^*(\alpha,r))^k = \alpha^k\cup\{y^{k-1}\}\setminus\{y^{k}\}
\end{equation}
\begin{figure}[htb]
\begin{center}
\input{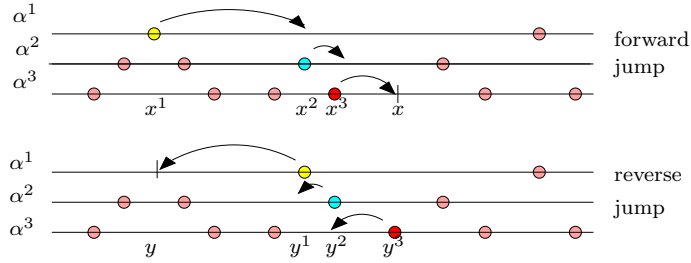}  
\label{330}
\caption{Effect of a point at $x$ for the multi-line process and the reverse jump
  (due to a point at $y$ in the reverse).}
\end{center}
\end{figure}

The operators $J$ and $J^*$ are the inverse of each other in the
following (equivalent) senses: 
\begin{equation}
  \label{f25}
  J^*(J(\alpha,x),x^1) = \alpha;\qquad J(J^*(\alpha,y),y^{n}) = \alpha
\end{equation}
(in the figure, $x^1=y$ and $y^n=x$.)  The generator of $\alpha^*_t$ is given by
\begin{equation}
  \label{f26}
  \cL^* f(\alpha) = \int_\R dy [f(J^*(\alpha,y))-f(\alpha)]
\end{equation}
\begin{propos}
  \label{f8}
  For any choice of positive densities $\rho=(\rho^1,\dots,\rho^n)$,
  the multi-line process with generator $\cL^*$ is the reverse of the
  multi-line process with generator $\cL$ with respect to the product
  measure $\nu$ on $\cX^n$ with densities $\rho$. As a consequence,
  $\nu$ is invariant for both $\cL$ and $\cL^*$.
\end{propos}

\paragraph{\bf Proof}
The proof is based on the observations: a) the configurations $\alpha$
and $J(\alpha,r)$, the configuration after a jump produced by a
Poisson point at $r$, have the same probability weight; b) the rate to
jump from $\alpha$ to $J(\alpha,r)$ is the same as the rate in the
reverse process to jump from $\alpha'=J(\alpha,r)$ to
$J^*(\alpha',y)=\alpha$; c) the rate of exiting any configuration
$\alpha$ is the same for the direct and the reverse process. We
formalize this with a generator computation.

It suffices to show that for bounded functions $f$ and $g$ depending
on finite regions in $\R^n$, $\nu (g \cL f)= \nu ( f \cL^*g)$; that
is, 
\begin{eqnarray}
  \label{f14}
  \int \nu(d\alpha)\int_\R g(\alpha) dx
  [f(J(\alpha,x)) - f(\alpha)]
  &=& \int \nu(d\alpha)\int_\R f(\alpha) dx
  [g(J^*(\alpha,x)) - g(\alpha)]\nonumber
\end{eqnarray}
Since the term subtracting in both terms is the same, it suffices to
prove that 
\begin{equation}
  \label{f33}
  \int \nu(d\alpha)\int_\R g(\alpha)f(J(\alpha,x)) dx =
\int \nu(d\alpha)\int_\R g(\alpha)f(J^*(\alpha,x)) dx 
\end{equation}
By conditioning on the positions of the jumping particles
$(x^1,\dots,x^{n})$ the left term in \reff{f33} can be written as
\begin{eqnarray}
 \label{f30}
 \int_\R  dx \,\int_{x^1<\dots< x^{n}<x}\, \prod_{k=1}^n
 (\rho^k e^{-\rho^k (x^{k+1}-x^k)})\, dx^1\dots dx^n \,
 \int \nu_{x}(d\alpha)\,g(\alpha)
 f(J^*(\alpha,x))
\end{eqnarray}
where $\nu_x$ is the measure $\nu$ conditioned to have at each line $k$
a particle at position $x^{k}$ and no particle in $(x^k,x^{k+1})$.
Change variables: call $\alpha' = J(\alpha,x)$ so that, by \reff{f25}
$\alpha= J^*(\alpha',y)$ and $y^k=x^{k+1}$, $k=0,\dots,n$ to get that
\reff{f30} equals
\begin{equation}
  \label{f34}
  \int_\R  dy \,\int_{y<y^1<\dots< y^{n}}
  \,\prod_{k=1}^n
 (\rho^k e^{-\rho^k (y^{k}-y^{k-1})})\, dy^1\dots dy^n \,
  \int \nu^{y}(d\alpha')\,g(J^*(\alpha',y))
  f(\alpha') 
\end{equation}
where $\nu^y$ is the measure $\nu$ conditioned to have at each line $k$ a
particle at position $y^k$ and no particles in $(y^{k-1},y^k)$.  Expression
\reff{f34} is just the right hand side of \reff{f33}.  \square

The set of \emph{dual points} of the \had{}  process starting with a Poisson process
of rate $\lambda$ is a Poisson process on $\R\times\R^+$.  This is proven in
Theorem 3.1 in \cite{CG} under the title ``Burke's theorem for Hammersley
processes''. We state and prove this fact in our context.
\begin{propos}
  \label{f11}
  If $\eta\in\cX$ has law $\nu^\lambda$ and $\omega$ is a rate-1 homogeneous
  Poisson process in $\R\times\R^+$, then the dual points of the trajectory
  $(\eta_t,\,t\ge 0)$ form a rate-1 homogeneous Poisson process in $\R\times\R^+$.
  Furthermore for $t>0$, $\eta_T$ is independent of the dual points contained in
  $\R\times[0,T]$.
\end{propos}

\paragraph{\bf Proof} 
By Proposition \ref{f8} (with $n=1$), the reverse \had{}  process with respect to
$\nu^\lambda$ is just a \had{}  process with the drift to the left. Assume we have
constructed the process $(\eta_t,\, 0\le t\le T)$ with initial law $\nu^\lambda$
and Poisson points $\omega$ on $\R\times[0,T]$.  The reverse process
$\eta^*_s=\eta_{T-s-}$, $0\le s\le T$, is governed by the dual points in
$\R\times [0,T]$. Since the reverse process started in equilibrium $\nu^\rho$,
the points governing the reverse process must be also Poisson in $\R\times
[0,T]$.

The dual points in $\R\times[0,T]$ are the points governing the future of the
reverse process $\eta^*_t:=\eta_{T-t-}$ which is in equilibrium.  Hence it is a
Poisson process independent of the initial configuration $\eta^*_{0}=\eta_T$.
\square

In the proof of Theorem \ref{t3} we use Proposition \ref{f11} to construct a
stationary trajectory $(\eta_t,\,t\in\R)$ of the \had{}  process governed by points
$\omega$ with marginal law $\nu^\lambda$. The dual points of this trajectory
are, as $\omega$, a homogeneous Poisson process in $\R^2$.

\section{Multi-line and multiclass processes}
\label{s3}
In this section we prove that the projection of the multi-line process on the
multiclass space using the map $M$ is just the multiclass process. The analogous
statement is true for the coupled process.  Recall that given an initial
multi-line configuration $\alpha$ and a homogeneous rate-1 Poisson process
$\omega$ in $\R\times\R^2$ the process $(\alpha_t,\,t\ge0)$, where
$\alpha_t=\psi(t,\alpha,\omega)$ is a realization of the multi-line process. On
the other hand, given an initial configuration $\xi$ and the same $\omega$, the
process $(\xi_t,\,t\ge0)$ with $\xi_t=\Phi(t,\xi,\omega)$ is the multiclass
process.

\begin{propos}
  \label{f18}
  Assume $C\alpha$ is well defined. Then almost surely:
\begin{eqnarray}
  \label{f20}
  C(\Psi(t,\alpha,\omega)) &=& \Phi(t,C\alpha,\omega)\\
  M(\Psi(t,\alpha,\omega)) &=& \Upsilon(t,M\alpha,\omega)\label{f201}
\end{eqnarray}
In particular, the law of $C\alpha_t$ is the same as the law of $\eta_t$.
\end{propos}

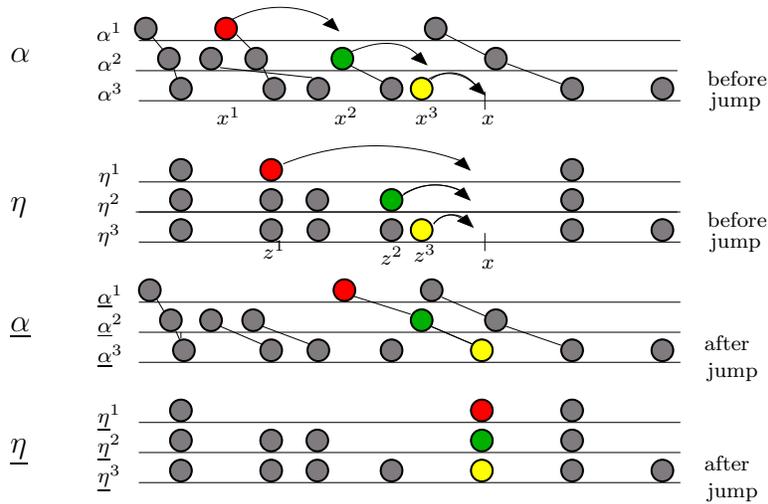
\begin{figure}[htb]
\begin{center}
\input{multiline-multitype.pstex_t}
\end{center}
\end{figure}
\paragraph{\bf Proof}
Expression \reff{f201} follows then from \reff{f20}, the identity $M=RC$ and
\reff{f6a}.  Fix $x\in\R$; it suffices to show
\begin{equation}
  \label{i0}
  H(C\alpha,x) = CJ(\alpha,x)
\end{equation}
where $H(\eta,x)$ is the coupled Hammersley configuration obtained after a point
at $(x,t)$ (for some $t$) if $\eta_{t-}=\eta$ and $J(\alpha,x)$ is the
resulting multi-line configuration after a point
at $(x,t)$ if $\alpha_{t-}=\alpha$: its $k$th coordinates satisfy
\begin{eqnarray}
  \label{f400}
  (H(\eta,x))^k &=& \eta  \setminus  \{z^k\} \cup \{x\} \nonumber\\
  (J(\alpha,x))^k &=& \alpha^k \setminus \{x^k\} \cup \{x^{k+1}\} \nonumber
\end{eqnarray}
where $z^k$ is position of the closest $\eta^k$ particle to the left of $x$ and
$x^k$ is the position of the particle jumping on line $k$ in the multi-line
process; here $x^{n+1}=x$.

Let $\eta=C\alpha$ be the coupled configuration obtained from $\alpha$ and
denote $\ualpha=J(\alpha,x)$ and $\ueta=H(\eta,x)$. With this notation \reff{i0}
is $C\ualpha=\ueta$. We first prove \reff{i0} for $n=2$ lines and then use this
case to proceed by induction.

\emph{The case $n=2$}. The second line of both processes are governed by the
point at $x$ with the same Hammersley rule. This implies they coincide in the
second line: $(C\ualpha)^2=(\ueta)^2$. There is an $\ualpha^2$ service time at
$x$ and an $\ualpha^1$ particle at $x^2$. This particle is served in line $2$ at
time $x$ unless there is another $\ualpha^1$ particle to the left of $x^2$
served at that time. In any case, an $\ualpha^1$ particle is served at $x$ in
line $2$. This means $(C\ualpha)^1$ has a particle at $x$. To see that
$(C\ualpha)^1$ has no particle at $z^1$ and that otherwise it is equal to $\eta$
we consider two cases depending whether $x^2$ serves or not an $\alpha^1$
particle.

1) The $\alpha^2$ time at $x^2$ serves an $\alpha^1$ particle. Hence, there is
an $\eta^1$ particle at $x^2$. Since there are no $\alpha^2$ services between
$x^2$ and $x$ (by definition of $x^2$), there cannot be $\eta^1$ particles
either and $\ueta^1= \eta^1\setminus\{x^2\}\cup\{x\}$.  Lets check that also
$C\ualpha$ is of this form. Let $y$ be the position of the $\alpha^1$ particle
served at $x^2$ in line 2. Clearly $y\le x^1$. The $\alpha^1$ particles to the
left of $y$ are served in line 2 before $x^2$ and hence its service time is not
affected if $x^1$ and $x^2$ are translated to the right. This means $C\ualpha$
coincides with $\ueta$ at $x$ and to the left of $x$. The $\alpha^1$ particles
to the right of $y$ are served after $x$. The displacement of the service time
from $x^2$ to $x$ just changes the service time of the particle at $y$, leaving
the other service times unchanged.  The displacement of the $\alpha^1$ particle
from $x^1$ to $x^2$ does not change the order or arrival of the particles
neither their service times at $\alpha^2$, as those times are after $x$. This
implies $C\ualpha=\ueta$ in this case.

2) The $\alpha^2$ time at $x^2$ does not serve an $\alpha^1$ particle. Then the
$\alpha^1$ particle at $x^1$ is served at some time $z$ before $x^2$. There are
no $\alpha^1$ particles between $x^1$ and $x^2$. The translation of the
$\alpha^1$ particle from $x^1$ to $x^2$ and the $\alpha^2$ particle from $x^2$
to $x$ do not change the service times of the $\alpha^1$ particles to the left
of $x^1$, which are served before $x^2$. By the FIFO schedule, particles arrived
after $x^2$ are served after particles arrived before. Hence, the translations
do not change the service times of the particles arrived after $x^2$.

\emph{The induction step}. From the definition of $C$, we have
\begin{equation}
\label{Tcorres}
\eta_t^k=D^{(n-k+1)}\left(\alpha^k_t,\dots,\alpha^n_t\right).
\end{equation}
From the observation just before \reff{Tdef}, $\eta_t^k$ has distribution
$\nu^{\rho^k}$.  So we need to show that the RHS of \reff{Tcorres} is a \had{} 
trajectory governed by $\omega$.  Since $(\alpha^k,\dots,\alpha^n)$ is itself a
multi-line process (with $n-k+1$ lines) governed by $\omega$, it is enough to
show that, for any $n$, $D^{(n)}(\alpha^1_t,\dots,\alpha^n_t)$ is a \had{} 
trajectory governed by $\omega$.  But from the definitions of $D^{(n)}$ and of
the multi-line process,
\[
D^{(n)}(\alpha^1_t,\dots,\alpha^n_t)=
D^{(2)}\left(
D^{(n-1)}(\alpha^1_t,\dots,\alpha^{n-1}_t), \alpha^n_t
\right)
\]
This and the fact that $(\alpha_t^1,\dots,\alpha_t^{n-1})$ is an $(n-1)$-line
multi-line process governed by $\omega^{n-1}$ concludes the induction step. \qed

\section{Uniqueness of the invariant measure}
\label{s5}
In this section we prove Theorems \ref{f5a} and \ref{t3} using the following
Proposition. It says that the marginals at time $t$ of the coupled Hammersley
process with initial independent marginals will coincide as $t\to\infty$ if the
first initial marginal is a Poisson process and the second one is an ergodic
distribution with the same density as the first.

\begin{propos}
  \label{t18}
  Let $\lambda>0$ and $\eta,\teta\in\cX^2$ be independent particle
  configurations with law $\nu^\lambda$, the Poisson process and
  $\tnu^\lambda$, an ergodic process in $\cX$ with density $\lambda$,
  respectively.  Let $(\eta_t,\teta_t)$ be the coupled process with
  initial configuration $(\eta,\teta)$. Then with probability one,
  both $\eta_t\setminus\teta_t$ and $\teta_t\setminus\eta_t$ have a
  deterministic density denoted $c(t)=c(t;\tnu^\lambda)$ which is decreasing and
  converges to 0 as $t\to\infty$. 
\end{propos}

An immediate consequence of the proposition and the translation invariance of
the law of the coupled process at each given time is that the probability of the
event $\{\Lambda\cap \eta_t= \Lambda\cap \teta_t\}$ converges to $1$ as
$t\to\infty$ for any bounded interval $\Lambda\subset \R$.

Before proving the proposition we
show how it implies Theorems \ref{f5a} and \ref{t3}.

\paragraph{\bf Proof of Theorem \ref{f5a}}
Assume $\eta$ and $\teta$ in $\cX^n$ are independent with laws $\mu$ and $\tmu$
respectively and let $(\eta_t,\teta_t)$ be a $2n$-coordinates coupled process
starting with configuration $(\eta,\teta)$. Assume $\mu$ and $\tmu$ are
invariant for the coupled process and that the marginal law of both $\eta^k$ and
$\teta^k$ is a Poisson process with rate $\rho^k$ for all $k$. Then, for bounded
$f$ depending on a bounded interval $\Lambda$:
\begin{eqnarray}
  \label{t19}
  |\mu f - \tmu f| &=& \Bigl|\int \mu(d\eta) f(\eta) - \int \tmu(d\teta)
  f(\teta)\Bigr|\\
  &=& \Bigl|\int \int \mu(d\eta)\tmu(d\teta)\, \E( f(\eta_t)-f(\teta_t))\Bigr|\\
  &\le& \int \int \mu(d\eta)\tmu(d\teta) ||f||_\infty \P(\Lambda\cap \eta_t\neq
  \Lambda\cap \teta_t)
\end{eqnarray}
which tends to zero by Proposition \ref{t18} applied to each
coordinate.  This implies $\mu=\tmu$.  \square

\paragraph{\bf Proof of Theorem \ref{t3}}
First construct a double infinite realization and the corresponding
points. Fix particles $\eta$ with law Poisson of rate $\lambda$ and
Poisson points $\omega^-\cup\omega^+$, the subsets of points with
negative and positive time coordinates, respectively. Run the process
$\eta_t=\Phi(t,\eta,\omega^+)$ for $t\ge 0$ using the points
$\omega^+$. Run the reverse process backwards $\eta^*_{-t} =
\Phi^*(-t,\eta,\omega^-)$ starting from the same configuration using
the points $\omega^-$. Here
$\Phi^*(-t,\eta,\omega^-)=\Phi(t,\eta,{\rm{TR}}(\omega^-))$, where
${\rm TR}(\omega)=\{(x,-t),\,(x,t)\in\omega\}$ are the points of
$\omega$ reflected with respect to the line $\{t=0\}$. Let the dual
points $D^-(\omega^-,\eta)$ be the positions of the particles of
$\eta^*_t$ just before jumps. By Proposition \ref{f11},
$D^-(\omega^-,\eta)$ is a Poisson process of points. Let
$\omega=D^-(\omega^-,\eta)\cup\omega^+$ be the configuration
consisting on the dual points of the reverse process for negative times
and the original points for positive times. The points $\omega$ are
Poisson and govern a stationary process having configuration $\eta$ at
time zero. This constructs simultaneously Poisson points $\omega$ and a
stationary trajectory governed by $\omega$.

For uniqueness we need to show that if
$(\teta_s)=(\teta_s,\,-\infty<s<\infty)$ is a stationary evolution
governed by $\omega$ such that the time marginal $\teta_s$ is ergodic
with density $\lambda$ then $\teta_s=\eta_s$ for all $s$. By Theorem
\ref{f5a} the marginal law of $\teta_t$ is Poisson of parameter
$\lambda$ for all $t$. Hence both $(\teta_s)$ and $(\eta_s)$ are
space-time ergodic processes but the joint process
$((\eta_s),(\teta_s))$ is not necessarily ergodic.  Fix a positive $t$
and introduce an auxiliary process $(\eta'_s,\,s\ge -t)$ governed by
$\omega$. The initial configuration $\eta'_{-t}$ is Poisson of rate
$\lambda$ and independent of the configurations
$(\eta_{-t},\teta_{-t})$ at that time.  By Proposition \ref{t18} the
density of $\eta_0\Delta\eta'_0$ and $\teta_0\Delta\eta'_0$ are both
smaller than $2c(t)$ for all $t$.  Hence the density of
$\teta_0\Delta\eta_0$ is bounded by $4c(t)$ for all $t$.  Since $c(t)$
converges to 0, this implies that $\eta_0=\teta_0$ and by the same
argument $\eta_s=\teta_s$ for all $s\in\R$.  \square

\paragraph{\bf Proof of Proposition \ref{t18}}
We use an argument of Ekhaus and Gray \cite{eg93} as developed in
Section 2 of Mountford and Prabhakar \cite{mp95}. Let
$\xi^=_t=\eta_t\cap\teta_t$, $\xi^+_t=\teta_t\setminus\eta_t$ and
$\xi^-_t=\eta_t\setminus\teta_t$.  As in \cite{mp95} we call them
\emph{yellow}, \emph{blue} and \emph{red} particles, respectively.
For each $t$ the process $(\xi^=_t,\xi^+_t,\xi^-_t)$ is ergodic and
$\xi^+_t$ has the same density as $\xi^-_t$. We want to show that the
density of $\xi^+_t$ goes to zero as $t\to\infty$.  Label the
particles of $\eta_t$ and $\teta_t$ as follows. Call $\eta_t(i)$ the
position of the $i$th particle of $\eta$. Initially
$\eta_0(i)<\eta_0(i+1)$ for all $i$; the same for $\teta$. The labels
evolve in time depending on the color. At the time of an $\omega$
Poisson point at $x$, proceed as follows:

1) If the closest left $\eta$ and $\teta$ particles are in $\xi^=$,
then both of them jump to $x$ carrying their labels.

2) If the closest left $\eta$ particle is blue localized in $x_1<x$,
the closest $\teta$ particle is red localized in $y_1<x_1$, then the
$\eta$ particle at $x_1$ and the $\teta$ particle at $y_1$ jump to $x$
carrying their labels and change their color to yellow. If furthermore
there are blue particles in $y_1<x_k<\dots<x_2<x_1$, the blue particle
in $x_i$ jumps to $x_{i+1}$, $i=2,\dots,k-1$ carrying the label and
the color.

2') The same as (2) by interchanging $\eta$ with $\teta$ and blue with
red.

3) If the closest $\eta$ particle is blue at position $x_1<x$ and the
closest $\teta$ particle is yellow at position $y_1<x_1$, then the
yellow particles in $y_1$ jump to $x$ carrying the label and keeping
the color. If furthermore there are blue particles in
$y_1<x_k,\dots,x_2<x_1$, then the $\eta$ blue particle at $x_i$ jumps
to $x_{i+1}$, $i=2,\dots,k-1$ carrying the label and the color.

3') The same as (3) by interchanging $\eta$ with $\teta$ and blue with
red.

In other words, when a labeled particle becomes yellow it keeps the
color for ever. Yellow particles behave as first class particles while
blue and red particles do as second class particles. When blue and red
particles coalesce they change the color to yellow; this happens in
cases (2) and (2'). Since there is no creation of new $\xi^{\pm}$
particles, its density must be non increasing. If a labeled $\eta$ or
$\teta$ particle overpasses another particle, then it becomes yellow.
As a consequence blue and red particles can be overpassed but cannot
overpass other particles.

Call \emph{everblue} those particles in the initial configuration $\eta$ that
will be blue at all $t\ge 0$ and \emph{everred} those particles in $\teta$ that
will be red at all times. The configuration of everblue particles has a
translation invariant distribution but not necessarily ergodic. Let
$\gamma^+_t\subset\eta$ be the set of $\eta$ particles that at time $t$ will be
blue. The law of $\gamma^+_t$ is ergodic with a (deterministic) density $c(t)$,
the density of $\xi^+_t$. $c(t)$ decreases to a value $c$ as $t\to\infty$.
Since the configuration of everblue particles is the intersection in $t$ of
$\gamma^+_t$, it has density $c$.  For the same reason the configuration of
everred $\teta$ particles has also ergodic law with density~$c$.

\begin{lemma}
  \label{b2}
  Assume the conditions of Proposition \ref{t18}. If the density of
  everblue particles is strictly positive, then there exist an
  $M<\infty$, an $n<\infty$ and a $\delta>0$ such that for each $t$,
  the density of red particles $W$ at time $t$ satisfying

\noindent {\rm i)} there exist blue particles at time $t$ in $(W,W+M]$,

\noindent {\rm ii)} there are at most $n$ $\eta_t$ particles in $(W,W+M]$,

\noindent {\rm iii)} there are no $\eta_t$ particles in $(W+M,W+M+\delta]$,

\noindent is at least  $1/(4M)$.
\end{lemma}

\paragraph{\bf Proof} The existence of an $M$ such that the the red
particles $W$ at time $t$ satisfying (i) have at least density $1/M$
is proven in Lemma 3.1 of \cite{mp95} for another process. This is the
hard part of the Ekhaus-Gray argument. The proof applies here because
it only uses the fact that red and blue particles cannot overpass, so
everred and everblue particles maintain the order.  Since $\eta_t$ is
a Poisson process, the density $a(k)$ of $\eta_t$ particles $U$ having at
least $k$ $\eta_t$ particles in the interval $(U,U+M)$ decreases
exponentially with $k$.  For $n$ sufficiently large so that
$a(n)<1/(2M)$ the density of red particles at time $t$ satisfying (i)
and (ii) for this $n$ is at least $1/(2M)$.  Finally, as $\eta_t$ is a
Poisson process, the density $b(\delta)$ of $\eta_t$ particles $U$
such that there is at least one $\eta$ particle in the interval $(U+M,
U+M+\delta]$ is bounded above by $\delta\rho$ (the mean number). Take
$\delta$ sufficiently small such that $b(\delta)<1/(4M)$. Then the
density of red particles $W$ at time $t$ satisfying (i), (ii) and such
that there is another $\eta$ particle in $(W+M,W+M+\delta]$ is smaller
than $1/(4M)$ and the red particles $W$ at time $t$ $\xi^-_t$ satisfying
(i-ii-iii) have at least density $1/(4M)$.  \square

\begin{lemma}
  \label{b3}
  Under the conditions of Lemma \ref{b2} there exists a positive
  $\vep'$ such that for all $t$ the density of $\xi^-_t$ minus the
  density of $\xi^-_{t+1}$ is at least $\vep'$.
\end{lemma}

\paragraph{\bf Proof} 
Take integer times $t$ and for each $\eta_t$ particle $U$ consider the
following event in the space-time Poisson process $\omega$:

\noindent iv) there are no $\omega$ points in
$(U,U+M]\times [t,t+1/2]$.

\noindent v) $\omega$ has exactly $n+1$ points
$(x_1,t_1),\dots,(x_{n+1},t_{n+1})$ in $(U+M,U+M+\delta]\times
[t,t+1]$ and they are increasing in the time coordinate and decreasing
in the space coordinate.

\noindent This event has a positive probability $\vep''$ and it is independent
of the past up to $t$. This implies that the density of red particles at time
$t$ satisfying conditions (i) to (v) has some positive density at least
$\vep'=\vep''/8M$.

The red particles at time $t$ satisfying (i) to (v) will collide with
a blue particle between $t$ and $t+1$. Hence the density of
$\xi^-_{t+1}-\xi^-_t$ is not smaller than $\vep'$. \square

\paragraph{\bf Proof of Proposition \ref{t18}}
By contradiction. If we assume that the density of $\xi^+_t$ decreases
to a non negative constant $c$, then by Lemma \ref{b3} this same
density decreases by a fixed amount $\vep'>0$ at each unit of time.
\square

\section{Regeneration properties of the multiclass invariant
  measure}
\label{s7}

Assume that there is a second class particle at the
origin. This corresponds to an unused service in the queue, which in
turn implies there are no customers in the queue at time $0$. Since
the queue is Markovian, the future depends only on the number of
customers at time 0. The attempted departures to the right of the
origin is a Poisson process with rate $\rho$. We conclude that
conditioned on having a second class particle at the origin, the first plus
second class particles to the right of the origin form a Poisson
process of rate $\rho$. By reversing time, the arrivals in the reverse
process also form a Poisson process of rate $\lambda$. The arrivals in
the reverse queue are just the effective departures in the forward
queue; that is, the first class particles. Hence, conditioned on having a
second class particle at the origin, the first class particles to the
left of the origin form a Poisson process of rate $\lambda$ and
independent of the right process.

For $n\ge3$ classes there are not regeneration events but as in the discrete
case, there are \emph{regeneration strings}. Fix a vector of classes
$(c_0,\dots,c_\ell)$ such that $c_0=n$, $c_\ell=2$ and for each class
$m\in\{3,\dots,n-1\}$ there exists a position $j_m<\ell$ with $c_{j_m}=m$ and
with $c_j\le m$ for all $j\in\{j_m+1,\dots,\ell\}$. For instance
$c=(4,1,2,3,1,2)$ qualifies for $\ell=5$ and $j_3=3$. We call $c$ a
\emph{regeneration string}. Let $\cX^n(c)$ be the set of configurations $\xi$ in
$\cX^n$ such that there are positions $0=x_0<\dots<x_\ell$ such that $x_i\in
\xi^{c_i}$ and $\xi^k\cap(x_i,x_{i+1})=\emptyset$ for $k=1,\dots,n$. In other
words, for $\xi\in\cX^n(c)$ there is a $\xi^n$ particle at the origin and
looking at the classes of the first $\ell$ particles to the right of the origin,
their classes follow the vector $c$. A translation invariant measure $\mu$ on
$\cX^n$ with densities $\rho^1,\dots,\rho^n$ conditioned on $\cX^n(c)$ is the
Palm measure $\hat\mu$ further conditioned on $\cX^n(c)$; here the Palm measure
is the measure conditioned on having a $\xi^n$ particle at the origin defined by
$\hat\mu f= (1/|I|\rho^n)\int d\mu(\xi)\sum_{r\in\xi^n\cap I} f(\theta_r\xi)$
for any measurable bounded set $I\subset \R$.  As in the discrete case one can
prove the following result
\begin{propos}
  \label{t10}
  Let $c$ be a regeneration string of length $\ell$ and
  $x_0<\dots<x_\ell$. Then, under the multiclass invariant measure $M\nu$
  conditioned on $\cX^n(c)$ the string $\xi\cap(-\infty,x_0)$ is independent of
  $\xi\cap(x_\ell,\infty)$.
\end{propos}

The proof is analogous to the two-classes case. One has to verify that after a
regeneration string the queue is empty, so whatever will happen in the future
(of $x_\ell$) depends on the Poisson processes of the future, so it is
independent of the past of $x_\ell$. In particular the superposition of the
classes to the right of the regeneration string
$\xi^1\cup\dots\cup\xi^n\cup(x_\ell,\infty)$ under the conditioned measure of
the Proposition is a Poisson process of rate $\rho^{n}$. See more details in
\cite{FM05}.

\section{Shocks in \had{} }
\label{s8}
Liggett \cite{Liggettcoupling} introduced the coupled process to prove ergodic
properties of the exclusion process; he called ``discrepancies'' what we call
second class particles. The interest on the study of invariant measures for the
coupled process was renewed when its association to shock measures become
clear. A shock related to a process $\eta_t$ is a (possibly random) position
depending on $t$ with the property that \emph{uniformly in time} the asymptotic
densities to the right and left of the shock are different. Wick \cite{Wick85}
showed the existence of that position for the totally asymmetric zero range
process with densities 0 and $\lambda$ to the left and right of the origin,
respectively.  In fact this is equivalent to the position of the leftmost
particle in the \tasep{} with densities 0 and $\lambda$ to the left and right of
the origin, respectively. Then \cite{FKS} and \cite{Fshock} (see also
\cite{Liggettbook2}) used the existence of an invariant measure for $\xi_t$ for
$n=2$ to show the existence of a shock in the \tasep. In a similar way Garcia
\cite{Garcia} proved the existence of a shock measure for the \had{} process.

We show now a way to construct a shock measure for the \had{} with asymptotic
densities $\lambda>\rho$ to the left and right of the origin, respectively.
Start with the invariant measure $\mu$ of Theorem \ref{f5} for the
coupled process with two marginals, $\nu^{\lambda}$ marginal for $\xi^1$ and $\nu^{\rho}$
marginal for $\xi^1\cup\xi^2$.  Let $\hat\mu$ be the measure $\mu$ conditioned to
have a second class particle at the origin, and $\xi$ be a two-class
configuration chosen according to $\hat\mu$. From $\xi$ construct a configuration
$\eta$ by superposing the negative particles of $\xi^1$ and all the particles of
$\xi^2$ except the one at the origin:
\begin{equation}
  \label{f50}
  \eta = [\xi^1\cap (-\infty,0)]\cup\xi^2\setminus\{0\}
\end{equation}
Call $S$ the map that transforms $\xi$ in $\eta$ and $\mu_{\lambda,\rho}$ the
law of $\eta$ so constructed. For the coupled process with initial
configurations $\eta$ and $\eta':=\eta\cup\{0\}$, there is only a discrepancy
for all times.  This discrepancy behaves like a second class particle. Call
$X_t$ its position.  The process $(\eta_t,X_t)$ is Markovian but the marginal
process $X_t$ is not. Let $\eta'_t$ be the process defined as the translation by
$X_t$ of $\eta_t$. Call $\phi'$ the operator that takes $t$, $\eta'$ and
$\omega$ in $\eta'_t$. Then $\mu'$ is invariant for $\eta'_t$. The proof is
based on the fact that $S$ commutes with the dynamics of $\xi_t$ in the
following sense:
\begin{equation}
  \label{t12}
  S \Phi^2(t,\omega,\xi) = \Phi'(t,\omega,S\xi) 
\end{equation}
where $\Phi^2$ is the operator that transforms a time $t$, points
$\omega$ and initial configuration $\xi$ in the configuration of the
two-class process at time $t$ as seen from the second class particle
$\tau_{X_t}\xi_t$. Here $\tau_x$ is the translation operator defined
by $\tau_x A = \{y-x,\,y\in A\}$. The map $\Phi'$ takes a time $t$,
points $\omega$ and an initial (one-class) configuration $\eta$ into
the process as seen from a second class particle $\tau_{X_t}\eta_t$.
This works in the same way as for the \tasep{}  \cite{FKS, Fshock} and for
the \had{}  \cite{Garcia}, so we omit the proof.

\section*{Acknowledgements}

Thanks to Eric Cator, Sheldon Goldstein and Herbert Spohn.

This paper is partially supported by FAPESP, CNPq, PRONEX. PAF thanks
hospitality and support from IHES and Laboratoire de Probabilit\'es of
Universit\'e de Paris 7. JBM thanks IME-USP, the Brazil-France
agreement and the FAPESP.

\parbox{0.52\linewidth}
{
\textsc{Instituto de Matem\'atica e Estat\'istica,\\
Universidade de S\~ao Paulo,\\
Caixa Postal 66281,\\
05311-970 S\~ao Paulo,\\
Brazil}\\
\texttt{pablo@ime.usp.br}\\
\texttt{http://www.ime.usp.br/$\tilde{\,\,\,\,\,}$pablo}}
\hfill
\parbox{0.44\linewidth}
{\textsc{Department of Statistics,\\
1 South Parks Road,\\
Oxford OX1 3TG,\\
United Kingdom}\\
\\
\texttt{martin@stats.ox.ac.uk}\\
\texttt{http://www.stats.ox.ac.uk/$\tilde{\,\,\,\,\,}$martin}}

\end{document}